\documentclass[aip,cha,reprint,twocolumn,groupedaddress,floatfix]{revtex4-2}
\usepackage{amsmath,amssymb}
\usepackage{tikz}
\usetikzlibrary{decorations.markings,arrows.meta}
\usepackage{graphicx}
\usepackage{hyperref}
\hypersetup{colorlinks=true, linkcolor=blue, citecolor=blue, urlcolor=blue}

\begin{document}

\title{Computable certificates for data-driven dimension reduction on
reaction-diffusion lattices}

\author{Yuhong Li}
\affiliation{School of Civil and Hydraulic Engineering, Huazhong University
of Science and Technology, Wuhan 430074, China}
\email{li_yuhong@hust.edu.cn}
\thanks{ORCID: 0009-0007-4329-4642.}

\begin{abstract}
Data-driven methods recover low-dimensional coordinates from
high-dimensional spatiotemporal data. Whether those coordinates track a
genuine invariant object or a finite-resolution artifact is not settled.
On reaction-diffusion lattices, the spectral-gap condition of
inertial-manifold theory reduces to a critical lattice size $n_*$.
Below $n_*$, the two-dimensional diffusion-map reading carries an existence
guarantee; above $n_*$ the guarantee lapses, and this lapse does not rule
out a slow manifold. The closed form for $n_*$ follows from the reaction
Lipschitz constant and the slow-mode growth rate. The certificate is tested
on a quadratic and a cubic model system. At the Brusselator boundary the
existence sign is checked by interval enclosure, and the spectral margin
decays as $n^{-4}$ across the transition. On strongly damped Hopf lattices
the slow manifold collapses to the uniform limit cycle; this degeneracy,
the trivial graph $\varphi\equiv0$, is reported directly rather than
concealed. The boundary gives a PDE-side prior for data-driven reduction:
it marks the resolution at which a recovered two-dimensional structure can
be trusted.
\end{abstract}


\maketitle

\noindent\textbf{High-dimensional simulations of coupled units often
conceal their slow behavior in a handful of coordinates. Reduction
algorithms can surface these coordinates automatically, yet it has
remained unclear whether they describe a real pattern of the dynamics or
a finite-resolution artifact. This work turns a classical condition from
the theory of low-dimensional attractors into a computable check for
reaction-diffusion lattices. From the reaction strength and the growth
rate of the slowest mode, one estimates a critical number of sites
$n_*$: at or below this size, a two-dimensional reading carries an
existence guarantee; above it, the guarantee lapses without proving that
no slow pattern exists. The boundary is tested on two standard nonlinear
systems, and the certificate margin decays predictably as the lattice
grows. The result offers a theory-informed prior: before trusting a
reduced model on a lattice, one can check whether the chosen resolution
falls inside the certified range.}

\section{Introduction}

Diffusion maps~\cite{Coifman2006,Berry2013},
dynamic mode decomposition~\cite{Williams2015,Korda2018a,Korda2018b}, and
deep-learning embeddings~\cite{Lusch2018} recover two- or
three-dimensional coordinates from simulations of thousands of coupled
elements. Reduced models for dissipative PDEs, built on approximate
inertial manifolds, pair this latent-variable discovery with
theoretical structure~\cite{Koronaki2024}. Whether those coordinates
describe a genuine slow manifold of the dynamics or a finite-resolution
artifact remains the open question.

Inertial-manifold theory answers this only conditionally.
Foias, Sell, and Temam~\cite{FST88} and Mallet-Paret and
Sell~\cite{MPS88} prove that under a spectral-gap condition on the
linearized operator, the long-time dynamics lie on a finite-dimensional
Lipschitz graph. Jolly~\cite{Jol89} constructs such a graph explicitly;
Bensoussan and Flandoli~\cite{Bensoussan1995} extend the theory to noise;
Robinson~\cite{Robinson01} gives the textbook treatment. The
condition stays abstract. It asks whether the slow-fast gap exceeds a
constant built from the nonlinear Lipschitz constant, but it does not
tell a practitioner simulating a 200-site lattice which side of the
threshold the lattice sits on.

Learning-side work approaches the same gap from the observable side.
Colbrook,
Mezi\'c, and Stepanenko~\cite{Colbrook2026}, extending
rigorous data-driven Koopman computation~\cite{ColbrookTownsend2024},
prove convergence and impossibility results for Koopman spectral
learning in abstract function spaces. Their analysis is system-agnostic:
it characterizes the geometry of the observable space but yields no
computable bound for a specific PDE discretization. Certified data-driven
models~\cite{Heide2026} bound learned low-order models by semidefinite
programming; periodic-orbit certification~\cite{Loh2026} proves
specific orbits near learned proxies. Neither asks at what resolution a
certified slow manifold exists on a concrete reaction-diffusion lattice.
Koopman invariance diagnostics~\cite{Conradie2026} quantify the closure
error of a learned dictionary, again on the data side; whether a
low-dimensional object exists at the simulated resolution remains open
on the PDE side.

The abstract condition can be evaluated on the lattice, with a
sufficient boundary at a critical size $n_*$. Below $n_*$ the
spectral-gap condition holds, the slow manifold is certified as a
Lipschitz graph, and diffusion maps recover a clean two-dimensional
structure. Above $n_*$ the condition fails and the graph loses its
guarantee; the lapse does not prove that no slow manifold exists. The
closed form for $n_*$ follows from the reaction Lipschitz constant and
the slow-mode growth rate. Interval enclosure at the Brusselator boundary
checks the existence sign, and the spectral margin decays as
$n^{-4}$, faster than the naive $n^{-2}$ spectral-gap scaling.

The certificate boundary recurs across nonlinear structure.
The Brusselator (quadratic
reaction, $n_*=60$) and the FitzHugh--Nagumo system (cubic reaction,
$n_*\approx 203$) differ by two orders of magnitude in the reaction
Lipschitz constant, which sets a factor-of-three difference in the
certified range. On strongly damped Hopf lattices the slow manifold
collapses to the uniform limit cycle itself. This degeneracy, the
trivial graph $\varphi\equiv 0$, is a structural feature rather than a
certificate failure: when diffusion dominates, no spatial structure
remains to learn. The certificate does not require the manifold to carry
nontrivial transverse structure. It still marks when the two-dimensional
object recovered by diffusion maps is the cycle and when the map reads
finite-resolution noise; on such lattices any two-dimensional output
tracks the cycle, not a nontrivial pattern.

The certificate supplies a PDE-side, theory-informed prior for reduction.
One computes $n_*$ from the reaction parameters first and checks whether
the simulation resolution lies inside the certified range, instead of
letting a learning algorithm find the slow manifold unaided. Inside it,
the two-dimensional coordinates carry the existence guarantee; outside
it, the certificate offers none.

Two systems test the certificate. The Brusselator, with $a=1$,
$D_1=D_2=1$, $h=0.02$, undergoes a Hopf bifurcation at $b=2$; at
$b=2.2$ the certificate gives $n_*\approx 60$. The FitzHugh--Nagumo
system near its Hopf point gives $n_*\approx 203$, its reaction
Lipschitz constant two orders of magnitude smaller. Diffusion maps on
simulated trajectories in both cases expose the predicted two
near-degenerate components inside $n_*$ and a degenerating spectrum
beyond it.

\section{Computable certificate}

Consider a periodic ring of $n$ sites with spacing $h$,
\begin{align}
\dot u_j &= f(u_j,v_j) + \frac{D_1}{h^2}(u_{j+1}-2u_j+u_{j-1}),\\
\dot v_j &= g(u_j,v_j) + \frac{D_2}{h^2}(v_{j+1}-2v_j+v_{j-1}).
\end{align}
The discrete Laplacian has eigenvalues
$\lambda_k = 4h^{-2}\sin^2(\pi k/n)$ for $k=0,\dots,n-1$. The $k=0$ mode
is spatially uniform and carries the slow dynamics (a Hopf pair or a
steady focus); every $k\ge 1$ mode is transverse.

Four quantities enter the certificate.

\noindent\textbf{(H1) Absorbing box.} A rectangle
$B=[u_{\min},u_{\max}]\times[v_{\min},v_{\max}]$ on which the reaction
field is forward invariant, obtained from the maximum-principle chain.
The box assumes a periodic ring with equal diffusion: a Neumann
boundary or unequal $D_1,D_2$ alters the maximum-principle chain, and
the invariant domain must be re-derived before the Lipschitz constant
is evaluated.

\noindent\textbf{(H2) Box Lipschitz constant.}
\begin{equation}
L_{\Box} = \sup_{(u,v)\in B}
\left\|\begin{pmatrix} f_u & f_v \\ g_u & g_v \end{pmatrix}\right\|_2 .
\end{equation}

This bound is global over $B$; a narrower invariant box, or a local
Lipschitz estimate over the reachable set, tightens $L_{\Box}$ and raises
$n_*$ if H1 still holds. Since $B$ is a convex rectangle, the segment
joining two points of $B$ lies in $B$, and the mean-value inequality gives
$\|F(x)-F(y)\|\le L_{\Box}\|x-y\|$ for $x,y\in B$; $L_{\Box}$ is thus a
global Lipschitz upper bound on the reaction field. The certificate uses the
Frobenius upper bound
$\|J\|_2\le\|J\|_F$, with Jacobian entries enclosed by interval
arithmetic; the global spectral-norm supremum is not machine-enclosed.

\noindent\textbf{(H3) Spectral gap.} The smallest transverse eigenvalue,
\begin{equation}
\mathsf{gap} = \lambda_1 = \frac{4}{h^2}\sin^2\!\left(\frac{\pi}{n}\right).
\end{equation}

On the uniform ring the transverse modes are the discrete Fourier
modes. A non-uniform mesh, a non-local diffusion, or an operator with a
weaker second transverse mode calls for recomputing the transverse
spectrum; the certificate then takes the actual smallest transverse
eigenvalue instead of $\lambda_1$ above.

\noindent\textbf{(H4) Slow-mode growth.}
The real part $\gamma$ of the linearized $k=0$ dynamics (positive for a
Hopf pair, negative for a stable focus). The certificate uses an upper
envelope of this real part: a larger $\gamma$ lowers the net transverse
damping $\mathsf{gap}/2-\gamma$, so the bound stays conservative and the
certified boundary is narrower rather than wider. On a stochastic lattice the
deterministic real part gives way to a mean-square or almost-sure slow
growth rate; the gap and Lipschitz inputs stay unchanged.

The Hadamard fixed-point argument for inertial manifolds requires the
transverse damping to dominate the nonlinear coupling. Define
\begin{equation}
C = \max\!\left\{\frac{2}{\mathsf{gap}},\;
\frac{1}{\mathsf{gap}/2 - \gamma}\right\},
\qquad
\kappa_{\Box} = C\, L_{\Box}.
\end{equation}

\paragraph{Certificate condition.}
The FST inertial-manifold existence theorem~\cite{FST88} applies under a
spectral-gap condition: the slow subspace attracts transverse directions
faster than the nonlinear coupling spreads across them. The theorem is not reproved here; its abstract condition is evaluated on
the lattice system (1)--(2) through the four inputs H1--H4.

\noindent\textbf{Certificate.} If the absorbing box $B$, the box Lipschitz
constant $L_{\Box}$, the transverse gap $\mathsf{gap}$, and the slow-mode
growth $\gamma$ satisfy
\begin{equation}
\frac{\mathsf{gap}}{2\kappa_{\Box}}>1,\qquad
\kappa_{\Box}=C L_{\Box},
\label{eq:cert}
\end{equation}
with
\begin{equation}
C = \max\!\left\{\frac{2}{\mathsf{gap}},\;
\frac{1}{\mathsf{gap}/2-\gamma}\right\},
\label{eq:C}
\end{equation}
then a two-dimensional Lipschitz slow manifold
$\mathcal{M}=\{(\xi,\varphi(\xi))\}$ exists on the lattice, with
$\mathrm{Lip}(\varphi)\le C$ and transverse attraction rate
$\mathsf{gap}/2-\gamma$. The condition is sufficient, not necessary: it
certifies existence for $n\le n_*$; a failure above $n_*$ means only
that the guarantee lapses, not that no slow manifold exists. The
components are evaluated below on a concrete lattice.

\paragraph{Derivation of the condition.}
The FST spectral-gap condition reduces to
\eqref{eq:cert} as follows; the computable inputs stay explicit. The
discrete Fourier modes form an orthogonal basis of
$X=\mathbb{R}^{2n}$. Decompose
$X=X_s\oplus X_\perp$ with $X_s$ the two-dimensional
$k=0$ subspace and $X_\perp$ its orthogonal complement of $k\ge1$
modes. The linear part block-diagonalizes as
$A=\mathrm{diag}(A_s,A_\perp)$ with
$\max\mathrm{Re}\,\mathrm{spec}(A_s)\le\gamma$ and
$\max\mathrm{Re}\,\mathrm{spec}(A_\perp)\le-\mathsf{gap}/2$. Here
$\mathsf{gap}=\lambda_1=4h^{-2}\sin^2(\pi/n)$ is the smallest
nonzero discrete-Laplacian eigenvalue; the transverse contribution of
the reaction Jacobian is not subtracted from this spectral bound but is
absorbed into the box Lipschitz constant $L_{\Box}$, so the damping term
is a conservative estimate and the net transverse attraction rate is
$\mathsf{gap}/2-\gamma$. The
nonlinear remainder satisfies $\|F(x)-F(y)\|\le L_{\Box}\|x-y\|$.

For two candidate graphs, let $q=Q(x_\varphi-x_\psi)$ and
$p=P(x_\varphi-x_\psi)$ be the fast and slow components of their
difference. Variation of constants on the fast equation gives
\begin{equation}
q(t)=e^{A_\perp t}q(0)+\int_0^t
e^{A_\perp(t-s)}Q\bigl[F(x_\varphi(s))-F(x_\psi(s))\bigr]\,ds .
\label{eq:varconst}
\end{equation}
The fast diffusion block is the symmetric discrete Laplacian, normal in
the 2-norm, so its semigroup norm equals the spectral bound and carries
no non-normal transient amplification;
$\|e^{A_\perp t}\|\le e^{-\mathsf{gap}t/2}$. The reaction difference
enters the variation-of-constants source term, bounded by
$L_{\Box}\|x_\varphi-x_\psi\|$; the slow component evolves under
$A_s$ with bound $e^{\gamma t}$. Two time integrals then fix the
contraction factors. The fast damping integrates to
$\int_0^\infty e^{-\mathsf{gap}t/2}\,dt=2/\mathsf{gap}$, so the coupling
the graph may carry is bounded by $2L_{\Box}/\mathsf{gap}$. Relative to
the slow growth $\gamma$, the fast directions decay at the net rate
$\mathsf{gap}/2-\gamma$, whose integral is
$1/(\mathsf{gap}/2-\gamma)$, giving the second factor
$L_{\Box}/(\mathsf{gap}/2-\gamma)$.

\begin{figure}[t]
\centering
\includegraphics[width=0.9\columnwidth]{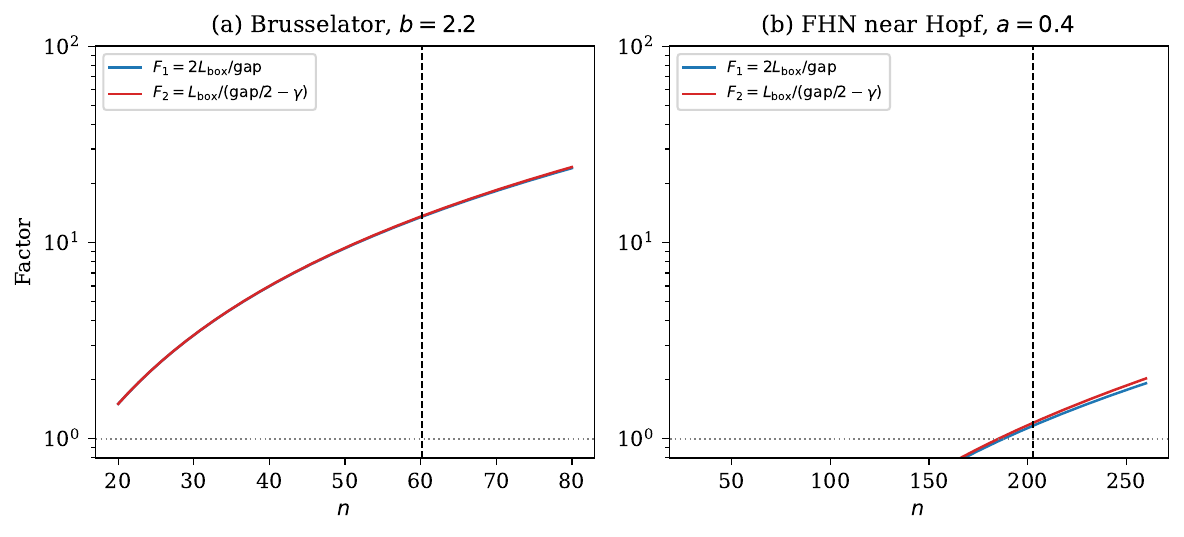}
\caption{The two contraction factors versus lattice size $n$.
(a) Brusselator $b=2.2$ ($L_{\Box}=184.36$, $\gamma=0.10$);
(b) FHN near Hopf $a=0.4$ ($L_{\Box}=1.40$, $\gamma=0.04$).
$F_1=2L_{\Box}/\mathsf{gap}$ is the graph-slope factor carried by the
transverse damping, and $F_2=L_{\Box}/(\mathsf{gap}/2-\gamma)$ is the
transverse-contraction factor relative to the slow growth. The dotted line
at factor $=1$ marks the stronger sufficient condition that each entry lie
below one; the certificate pass line is the ratio
$\mathsf{gap}/(2\kappa_{\Box})>1$, marked by the vertical dashed line at
$n_*$.}
\label{fig:factors}
\end{figure}

These are the two competing contraction factors
(Figure~\ref{fig:factors}). The first,
$2L_{\Box}/\mathsf{gap}$, measures how strongly the transverse damping
$\mathsf{gap}$ suppresses the nonlinear coupling $L_{\Box}$ when the
graph is flattened onto the slow subspace. The second,
$L_{\Box}/(\mathsf{gap}/2-\gamma)$, measures how fast transverse
directions contract relative to the slow growth $\gamma$; this factor
grows whenever $\gamma$ approaches $\mathsf{gap}/2$. The Hadamard
fixed-point argument constructs $\varphi$ as the fixed
point of a map $\mathcal{T}$ on Lipschitz graphs. The contraction
constant of $\mathcal{T}$ is bounded by
$(2L_{\Box}/\mathsf{gap})\cdot(L_{\Box}/(\mathsf{gap}/2-\gamma))$.
Requiring this product below $1$ is the contraction condition. To express
it through a single graph Lipschitz constant, the two factors are bounded
separately and each required below $1$; this yields
$C=\max\{2/\mathsf{gap},\,1/(\mathsf{gap}/2-\gamma)\}$ and the certificate
condition $\mathsf{gap}/(2C L_{\Box})>1$. The separate-entry requirement
is a conservative strengthening of the product condition: each factor
below $1$ suffices for the product to be below $1$, but the converse need
not hold, so the certified boundary is correspondingly narrower rather
than wider. The two
entries of $C$ track the two factors separately: $2/\mathsf{gap}$
bounds the graph slope allowed by the transverse damping, and
$1/(\mathsf{gap}/2-\gamma)$ bounds it by the transverse contraction
relative to the slow growth. This is the lattice form of the standard FST
spectral-gap condition, evaluated explicitly on the lattice.

\medskip
The condition is explicit in $n$. The certificate inequality
$\mathsf{gap}>2C L_{\Box}$ splits over the two entries of $C$. The first
entry $C_1=2/\mathsf{gap}$ gives
$\mathsf{gap}>2(2/\mathsf{gap})L_{\Box}$, hence
$\mathsf{gap}^2>4L_{\Box}$ and $\mathsf{gap}>2\sqrt{L_{\Box}}$. The
second entry $C_2=1/(\mathsf{gap}/2-\gamma)$, for a Hopf slow mode with
$\gamma>0$, gives
$\mathsf{gap}>2L_{\Box}/(\mathsf{gap}/2-\gamma)$, hence
$\mathsf{gap}^2-2\gamma\,\mathsf{gap}-4L_{\Box}>0$, whose positive root
is $\mathsf{gap}>\gamma+\sqrt{\gamma^2+4L_{\Box}}$. For $\gamma>0$ this
root is the tighter of the two thresholds. Substituting the large-$n$
asymptotics $\mathsf{gap}\approx 4\pi^2/(h^2 n^2)$ gives the approximate
critical lattice size
\begin{equation}
n_*^{\rm approx} = \frac{2\pi}{h\sqrt{\gamma+\sqrt{\gamma^2+4L_{\Box}}}}.
\label{eq:nstar}
\end{equation}
The expression is asymptotic. Since
$\sin x=x-x^3/6+O(x^5)$,
$\sin^2(\pi/n)=(\pi/n)^2[1-\pi^2/(3n^2)+O(n^{-4})]$; replacing the sine
by its leading term perturbs the gap by a relative $O(n^{-2})$, so the
closed form holds when $n$ is large enough that $\pi/n\ll1$. The exact
boundary
solves $\mathsf{gap}(n)/(2\kappa_{\Box}(n))=1$ numerically,
with $\mathsf{gap}(n)=4h^{-2}\sin^2(\pi/n)$ the exact discrete
Laplacian eigenvalue. The closed form guides intuition; the numerical
solve gives the certified boundary. The two agree within
1\% at the Brusselator operating point ($n_*^{\rm approx}=60.2$,
numerical boundary $n\approx 60$).

\paragraph{Scaling.}
Eq.~\eqref{eq:nstar} gives
$n_*\propto h^{-1}L_{\Box}^{-1/4}$ in the strongly damped regime
($|\gamma|\ll\sqrt{L_{\Box}}$): the inner square root scales as
$L_{\Box}^{1/2}$ and the outer one as $L_{\Box}^{1/4}$. The certificate ratio scales as
$\mathsf{gap}/(2\kappa_{\Box})\propto n^{-4}$: the extra $n^{-2}$ beyond
the naive $n^{-2}$ spectral-gap scaling comes from the factor
$C\sim 2/\mathsf{gap}\propto n^2$ in the coupling constant. A log-log
fit to the five certified points gives a slope of $-3.99\pm 0.08$
($R^2=0.999$), consistent with the predicted $-4$.

\paragraph{Proposition (trivial graph).}
When $D_1=D_2$ and the reaction vector field is invariant under
spatial translation (so that the uniform steady state or uniform limit
cycle is an exact solution), the transverse projection of that uniform
orbit is identically zero. The inertial graph coincides with the
embedding of the slow orbit into the full phase space, and
$\varphi(\xi)\equiv 0$.

\textit{Proof.} Let $\bar x(t)$ be the uniform orbit, i.e.\
$\bar u_j(t)=\bar u(t)$ and $\bar v_j(t)=\bar v(t)$ for all $j$. The
discrete Laplacian vanishes on $\bar x$ because all neighbors are
equal. The full lattice state $x(t)=(\bar x(t),\dots,\bar x(t))$
satisfies the lattice ODE exactly. Its projection onto the
$k\ge1$ transverse subspace is zero by construction. Since the FST
theorem guarantees a unique Lipschitz slow manifold in the
certified regime, this exact solution must lie on that manifold.
The graph $\varphi$ therefore satisfies $\varphi(\bar x(t))=0$ for
all $t$, and by Lipschitz continuity and uniqueness
$\varphi(\xi)\equiv 0$ on the slow subspace. \hfill$\square$

The Brusselator at $b=2.2$ sits in this case: the Hopf limit cycle
is spatially uniform, and no transverse mode is spontaneously excited.
The certificate still applies. It guarantees that the dynamics create
no transverse component, while the manifold it certifies is the cycle
itself.

\paragraph{Corollary (steady-state slow mode).}
When the $k=0$ mode is a stable focus ($\gamma<0$ rather than a Hopf
pair), the slow subspace is one-dimensional rather than two-dimensional,
and the attraction rate is $\mathsf{gap}/2+|\gamma|$ instead of
$\mathsf{gap}/2-\gamma$. The second entry of $C$ becomes
$1/(\mathsf{gap}/2+|\gamma|)$. The certificate inequality over this
entry gives
$\mathsf{gap}^2+2|\gamma|\,\mathsf{gap}-4L_{\Box}>0$, whose positive
root is $\mathsf{gap}> -|\gamma|+\sqrt{\gamma^2+4L_{\Box}}$. The first
entry still imposes $\mathsf{gap}>2\sqrt{L_{\Box}}$. With $g=|\gamma|$,
write the second threshold $t_2=-g+\sqrt{g^2+4L_{\Box}}$; the inequality
$\sqrt{g^2+4L_{\Box}}\le g+2\sqrt{L_{\Box}}$ reduces on squaring to
$0\le 4g\sqrt{L_{\Box}}$, so $t_2\le 2\sqrt{L_{\Box}}$. The first entry
is therefore the tighter one on the stable-focus branch, and the closed-form
boundary is
\begin{equation}
n_*^{\rm ss} = \frac{2\pi}{h\sqrt{2\sqrt{L_{\Box}}}}
=\frac{2\pi}{h(4L_{\Box})^{1/4}}.
\label{eq:nstar_ss}
\end{equation}
The second-entry expression $-|\gamma|+\sqrt{\gamma^2+4L_{\Box}}$ is the
looser slow-growth estimate; it under-reads the threshold and over-reads
$n_*$. At fixed $L_{\Box}$ the stable-focus threshold $2\sqrt{L_{\Box}}$
lies below the Hopf threshold $\gamma+\sqrt{\gamma^2+4L_{\Box}}$, so
$n_*^{\rm ss}>n_*^{\rm approx}$: a stable steady state carries no slow
growth rate competing with the transverse damping.
\hfill$\square$

\paragraph{Remark (noise robustness).}
The certificate is deterministic: it bounds manifold existence in the
noiseless lattice dynamics. Observation noise of relative amplitude
$\varepsilon_{\rm obs}$ perturbs the diffusion-map embedding but leaves
the manifold itself untouched; the manifold belongs to the dynamics,
not the data. While the noise level stays small enough for the embedding
bandwidth $\sigma$ to separate the near-degenerate slow pair from the
transverse modes, the two-dimensional reading holds. Sec.~IV checks
this numerically at $n=30$ across $\varepsilon_{\rm obs}\in\{0.5\%,1\%,5\%\}$.

\paragraph{Remark (Turing instability vs.\ certificate failure).}
The certificate fails for two distinct reasons. Resolution failure:
$n>n_*$ leaves the transverse damping $\lambda_1$ too small relative to
the nonlinear coupling, and the spectral-gap excess
$\mathsf{gap}/(2\kappa_{\Box})-1$ turns negative. Dynamical failure: a
transverse eigenvalue crosses the imaginary axis
($\mathrm{Re}\,\lambda_k>0$ for some $k\ge1$), a Turing instability
rather than a resolution failure. The two are distinguishable in the
certificate: a Turing instability invalidates the absorbing set
(the orbit leaves $B$); a resolution failure keeps the orbit inside $B$
but violates the contraction condition. The Brusselator at $b=2.2$
fails on resolution, not Turing.

Eq.~\eqref{eq:nstar} states the computable certificate: given $h$,
$\gamma$, and $L_{\Box}$, lattice sizes $n\le n_*$ are certified. The
bound is conservative; the actual certified-to-uncertified transition
sits near $n_*$ and is checked numerically below.

The conservatism enters at three separate steps. The product
condition $(2L_{\Box}/\mathsf{gap})(L_{\Box}/(\mathsf{gap}/2-\gamma))<1$ is
replaced by requiring each factor below one; each factor below one suffices
for the product to fall below one, but the converse need not hold, so the
certified set is a subset of what the product condition would allow. The
fast spectral rate $-\mathsf{gap}/2$ does not subtract the
transverse contribution of the reaction Jacobian; that contribution is
absorbed wholesale into $L_{\Box}$, which raises the Lipschitz constant and
weakens the damping term. The box Lipschitz constant and absorbing box are
analytic upper bounds---the interval-enclosed Frobenius norm estimate and
the maximum-principle chain---rather than machine-optimized global extrema. Each
step pushes $n_*$ downward; none widens the certified range. A failure above
$n_*$ means the guarantee lapses, not that no slow manifold exists.

Certificate readings use outward interval arithmetic with
60-digit endpoints, not single-point floating values.
At the tightest Brusselator certified point $n=60$, $b=2.2$, the gap ratio
lies in $[1.00995,1.00995]$, whose lower bound exceeds $1$ by $0.010$; the
FHN point $n=200$ gives $[1.0517,1.0517]$, and the uncertified points
enclose below $1$. The existence sign therefore does not hinge on
floating-point roundoff. This interval check covers the roundoff layer
only. The box Lipschitz constant and the absorbing box remain analytic
upper bounds (the corner norm estimate and the maximum-principle
chain), and are not verified by machine-assisted global optimization; the
readings stay a computable sufficient condition evaluated numerically,
not a strict computer-assisted proof of the Lipschitz constants or the
box.

\paragraph{Worked example (Brusselator, $b=2.2$).}
The certificate inputs step through as follows. The absorbing box
uses $m=a/(b+1)=1/3.2=0.3125$, $u_{\max}=a+b/m=1+2.2/0.3125=8.04$,
$v_{\max}=b/m=7.04$. The box Lipschitz constant is
$L_{\Box}=\sqrt{2}\,u_{\max}\sqrt{u_{\max}^2+4v_{\max}^2}=184.4$. The
slow growth is $\gamma=(b-2)/2=0.1$. At $n=60$,
$\mathsf{gap}=4/h^2\sin^2(\pi/60)=10000\cdot\sin^2(3^\circ)=27.39$.
The coupling constant is
$C=\max\{2/27.39,\,1/(27.39/2-0.1)\}=\max\{0.0730,\,0.0736\}=0.0736$.
The certificate ratio is
$\mathsf{gap}/(2 C L_{\Box})=27.39/(2\cdot0.0736\cdot184.4)=1.010$.
The lattice passes by a margin of 1\%. At $n=61$ the ratio drops below
1, so the numerical boundary is $n\approx 60$; the closed-form estimate
Eq.~\eqref{eq:nstar} gives $n_*^{\rm approx}=60.2$.

\section{Two model systems}

Table~\ref{tab:params} lists the parameters for both systems.

\begin{table}[t]
\caption{Model parameters. Both systems use $D_1=D_2=1$, $h=0.02$,
periodic boundary conditions on a ring. For FHN, $\gamma=0.04$ is a
conservative upper bound on the linearized real part $0.036$.}
\label{tab:params}
\begin{tabular}{l r r}
\small & \small Brusselator & \small FHN\\
\hline
\small Reaction & \small $a=1$, $b=2.2$ & \small $a=0.4,b=0.8,\varepsilon=0.08$\\
\small Hopf location & \small $b=2$ & \small $1-u_*^2-\varepsilon b=0$\\
\small Slow growth $\gamma$ & \small $0.1$ ($b=2.2$) & \small $0.04$ (cons.\ bound)\\
\small Box $L_{\Box}$ & \small $184.4$ & \small $1.4$\\
\small Certified $n_*$ & \small $\approx 60$ & \small $\approx 203$\\
\small Observation noise & \small $1\%$ rel. & \small $1\%$ rel.\\
\end{tabular}
\end{table}

\subsection{Brusselator}

The Brusselator reaction is
\begin{equation}
f(u,v)=a-(b+1)u+u^2 v,\qquad g(u,v)=bu-u^2 v,
\end{equation}
with $a=1$, $D_1=D_2=1$, $h=0.02$. The steady state $(a,b/a)$ undergoes a
Hopf bifurcation at $b=2$, with slow-mode growth $\gamma=(b-2)/2$. The
maximum-principle box uses $m=a/(b+1)$, $u_{\max}=a+b/m$,
$v_{\max}=b/m$. The reaction Jacobian is not rank one in general: its
determinant is $u^2(b+1)>0$ inside the box. The Frobenius norm, which
bounds the 2-norm, is evaluated at the corner $(u_{\max},v_{\max})$,
\begin{equation}
L_{\Box}=\sqrt{2}\,u_{\max}\sqrt{u_{\max}^2+4v_{\max}^2}.
\end{equation}
At $b=2.2$, $L_{\Box}=184.4$ and Eq.~\eqref{eq:nstar} gives
$n_*^{\rm approx}=60.2$ (numerical boundary $n\approx 60$).

Table~\ref{tab:nstar_b} gives the boundary across the
Hopf neighborhood $b\in[2.05,2.40]$. Moving away from the Hopf point
pulls $n_*$ inward from $63.6$ to $56.1$. Two inputs shift together:
the slow growth $\gamma=(b-2)/2$ rises against the transverse
damping, and the absorbing box widens with $b$, raising
$L_{\Box}$. On the Brusselator the box effect dominates; the
reaction coupling nearly doubles over the range, and the boundary
tracks the box Lipschitz constant rather than the growth rate alone.

\begin{table}[t]
\caption{Certificate boundary $n_*^{\rm approx}$ across the Hopf
neighborhood for the Brusselator. $\gamma=(b-2)/2$, $L_{\Box}$ is the
box Lipschitz constant.}
\label{tab:nstar_b}
\begin{tabular}{r r r r}
$b$ & $\gamma$ & $L_{\Box}$ & $n_*^{\rm approx}$\\
\hline
2.05 & 0.025 & 148.3 & 63.6\\
2.10 & 0.050 & 159.6 & 62.4\\
2.20 & 0.100 & 184.4 & 60.2\\
2.30 & 0.150 & 211.9 & 58.1\\
2.40 & 0.200 & 242.4 & 56.1
\end{tabular}
\end{table}

\subsection{FitzHugh--Nagumo}

The FitzHugh--Nagumo reaction is
\begin{equation}
f(u,v)=u-u^3/3-v,\qquad g(u,v)=\varepsilon(u+a-bv),
\end{equation}
with $D_1=D_2=1$, $h=0.02$. The steady state $(u_*,v_*)$ solves
$u_*-u_*^3/3-v_*=0$ and $u_*+a-bv_*=0$. With $a=0.4$, $b=0.8$,
$\varepsilon=0.08$, the steady state is $u_*\approx-0.93$,
$v_*\approx-0.66$. The reaction Jacobian there is
\begin{equation}
J = \begin{pmatrix} 1-u_*^2 & -1 \\ \varepsilon & -\varepsilon b \end{pmatrix}
\approx \begin{pmatrix} 0.136 & -1 \\ 0.08 & -0.064 \end{pmatrix}.
\end{equation}
The trace is $0.072$ and the determinant is $0.071$, giving a complex
pair with real part $0.036$. The Hopf threshold is where the trace
crosses zero, i.e.\ where $1-u_*^2-\varepsilon b=0$; the operating
point sits just on the unstable side, so the slow growth has linearized
real part $\gamma_{\rm lin}\approx 0.036$ (half the trace, for a
$2\times2$ Jacobian with complex eigenvalues). The certificate uses the
conservative upper bound $\gamma=0.04$. The box Lipschitz constant is
$L_{\Box}\approx 1.4$, carried by the Frobenius upper bound on the
reaction Jacobian. Equation~\eqref{eq:nstar} gives
$n_*\approx 203$, agreeing with the numerical boundary.

The same cubic reaction reaches a stable-focus operating point. With
$a=0.5$, $b=0.8$, $\varepsilon=0.08$, the steady state moves to
$u_*\approx-1.03$, $v_*\approx-0.67$, and the reaction Jacobian is
\begin{equation}
J \approx \begin{pmatrix} -0.066 & -1 \\ 0.08 & -0.064 \end{pmatrix},
\end{equation}
with trace $-0.130$, determinant $0.084$, and a complex pair of real part
$\gamma_{\rm lin}\approx-0.065$. The certificate uses the conservative bound
$\gamma=-0.06$; the box Lipschitz constant is
$L_{\Box}\approx1.22$, carried by the Frobenius upper bound.
Equation~\eqref{eq:nstar_ss} gives
$n_*^{\rm ss}\approx212$, matching the ratio scan, which crosses $1$ near
$n\approx212$.
The certificate passes at $n=200$ (ratio $1.25$) and lapses by $n=240$
(ratio $0.60$), with the existence sign enclosed by outward interval
arithmetic. This verifies the steady-state branch on a focus rather than a
Hopf pair.

The systems differ in nonlinear structure (quadratic versus cubic) and
in $L_{\Box}$ (two orders of magnitude), which accounts for the
factor-of-three difference in $n_*$. The comparison shows the certificate
is model-independent: the certified range follows the reaction Lipschitz
constant, not the particular nonlinearity.

\section{Data-driven recovery on certified points}

Each lattice is simulated with a semi-implicit Fourier scheme; 1\%
relative observation noise is added, and diffusion maps (Gaussian kernel,
bandwidth $\sigma=1.2$, six retained coordinates) run on a subsample of 3000
states from 12 training trajectories. The leading DM eigenvalues mark
the slow geometry empirically: a two-dimensional slow manifold reads as
two near-degenerate leading eigenvalues followed by a gap.

\subsection{Brusselator across the certified range}

Table~\ref{tab:dm} shows the leading diffusion-map eigenvalues at $b=2.2$
for $n=20$ through $n=60$, alongside the certificate reading
$\mathsf{gap}/(2\kappa_{\Box})$.

\begin{table}[t]
\caption{Brusselator at $b=2.2$: leading diffusion-map eigenvalues and
certificate ratio across the certified range. The first DM eigenvalue is
normalized to one and omitted.}
\label{tab:dm}
\begin{tabular}{r r r r r}
$n$ & $\lambda_1^{\rm DM}$ & $\lambda_2^{\rm DM}$ &
$\lambda_3^{\rm DM}$ & gap ratio\\
\hline
20 & 0.962 & 0.934 & 0.858 & 81.1\\
30 & 0.973 & 0.964 & 0.900 & 16.2\\
40 & 0.979 & 0.971 & 0.923 & 5.12\\
50 & 0.983 & 0.979 & 0.938 & 2.10\\
60 & 0.985 & 0.970 & 0.949 & 1.01
\end{tabular}
\end{table}

The two leading eigenvalues stay near-degenerate
at every certified point, consistent with the two-dimensional Hopf slow
manifold. The gap between $\lambda_2^{\rm DM}$ and
$\lambda_3^{\rm DM}$ shrinks monotonically from $0.076$ at $n=20$ to
$0.021$ at $n=60$, tracking the certificate margin
$\mathsf{gap}/(2\kappa_{\Box})$ as it drops from $81$ to $1.01$. The
empirical spectrum and the theoretical certificate deteriorate together.

Across the five certified points, the correlation between the certificate
margin
$\log(\mathsf{gap}/(2\kappa_{\Box}))$ and the DM spectral gap
$\lambda_2^{\rm DM}-\lambda_3^{\rm DM}$
is $r=0.97$ (Pearson). The correspondence is not
one-to-one. At $n=80$, where the certificate margin drops to $0.32$,
the DM gap still reads $0.029$, larger than at the certified $n=60$ point
($0.021$); the DM spectrum deteriorates more slowly than
the certificate predicts. The certificate is sufficient, not necessary,
so a lattice can stay data-recoverable after the guarantee lapses.

\begin{figure*}[t]
\centering
\includegraphics[width=0.95\textwidth]{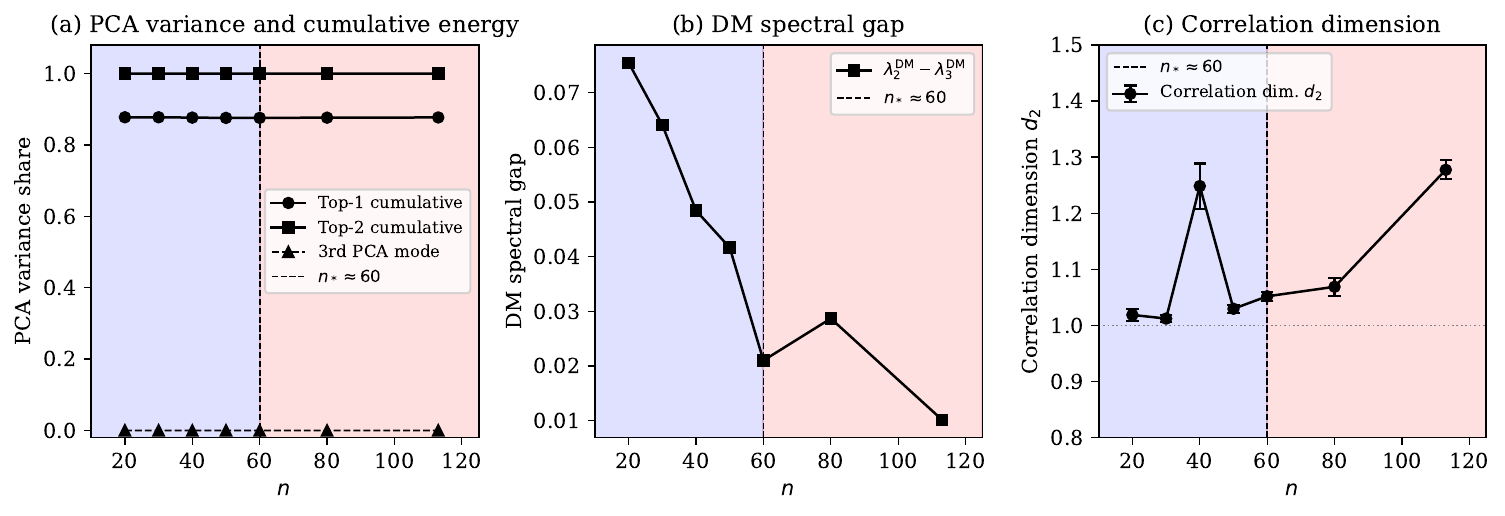}
\caption{Baseline comparison across the certified and diagnostic arms.
(a) PCA variance and cumulative energy: the top-1 and top-2 cumulative
shares stay flat on both sides of $n_*$, and the third mode stays
$\approx 0$; the energy criterion gives no sharp truncation and does not
detect the boundary. (b) DM spectral gap
$\lambda_2^{\rm DM}-\lambda_3^{\rm DM}$, which shrinks toward the
boundary. (c) Correlation dimension $d_2$ from the correlation integral
in the two leading PCA coordinates (three seeds, error bars $\pm 1$
standard deviation): it fluctuates about the cycle dimension $1$ on both
sides of $n_*$, with scatter already inside the certified arm, and does
not track the boundary. The PCA variance and correlation dimension
describe the observed trajectory; the DM gap responds to the certificate
margin.}
\label{fig:pca_dm}
\end{figure*}

PCA serves as a baseline against diffusion maps. On the certified
$n=30$ point, the leading PCA variance shares are
$0.877$, $0.123$, and $\approx 0$ from the third mode onward: PCA
recovers a two-dimensional structure, and the drop after the second
mode ($0.123$ to $\approx 0$) is sharper than the DM drop
($0.964$ to $0.900$). On the diagnostic arm the PCA shares remain
$0.876/0.124$ at $n=80$ and $0.877/0.123$ at $n=113$, with the third
mode still $\approx 0$. The variance concentration sits on the
two uniform modes of the limit cycle and does not detect the boundary at
$n_*$. A third intrinsic-dimension estimate, the correlation dimension
from the correlation integral in the two leading PCA coordinates, sits
near the cycle dimension $1$ on both arms ($d_2\approx1.0$--$1.3$ over
three seeds) and likewise fails to track $n_*$; the local value at
$n=40$ ($1.25\pm0.04$) falls inside the seed-to-seed scatter rather than
reading as a monotone boundary signal, and that scatter appears
already inside the certified arm. Only the DM spectral gap responds to the
certificate margin; it shrinks on the diagnostic arm ($0.029$ at $n=80$,
$0.010$ at $n=113$). PCA variance describes the observed trajectory,
whereas the certificate bounds manifold existence on the lattice;
the two should not be conflated. Figure~\ref{fig:pca_dm} contrasts
the two baselines across the certified and diagnostic arms.

\subsection{Robustness}

The analysis at $n=30$ is perturbed across observation noise (0.5\%, 1\%,
5\%), diffusion-map bandwidth ($\sigma=0.8,1.2,2.0$), and three trajectory
seeds. The $\lambda_2^{\rm DM}-\lambda_3^{\rm DM}$ gap stays positive in
all nine settings, ranging from $0.034$ (narrowest, at $\sigma=0.8$) to
$0.064$. The two-dimensional reading holds in every setting.

\begin{figure}[t]
\centering
\includegraphics[width=0.8\columnwidth]{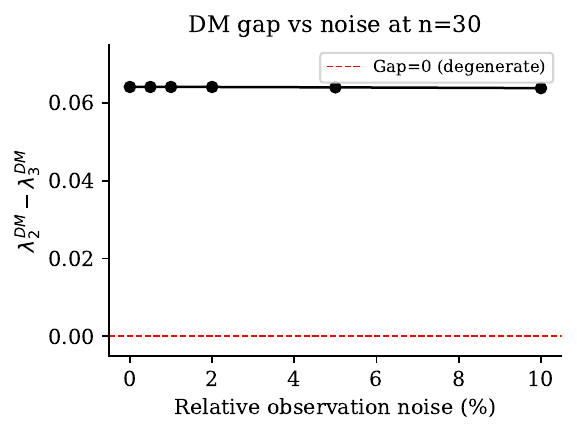}
\caption{The diffusion-map spectral gap
$\lambda_2^{\rm DM}-\lambda_3^{\rm DM}$ as a function of relative
observation noise at $n=30$ (certified). The gap stays near $0.064$
and remains positive out to $10\%$ noise; the two-dimensional reading
is stable.}
\label{fig:noise}
\end{figure}

Figure~\ref{fig:noise} plots the DM spectral gap against
observation noise at $n=30$ over a $0$--$10\%$ relative-noise
scan. The gap stays flat, $0.0641$ at zero noise to
$0.0638$ at $10\%$, and remains positive throughout; the
two-dimensional structure survives the tested range.

\begin{figure}[t]
\centering
\includegraphics[width=0.8\columnwidth]{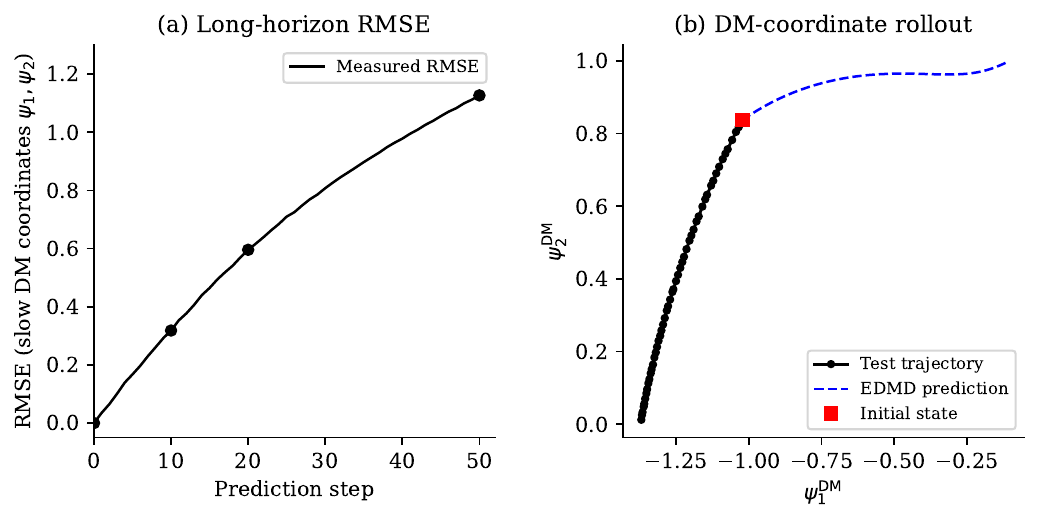}
\caption{EDMD long-horizon rollout at the certified $n=30$ point.
(a) RMSE on the two slow DM coordinates grows approximately linearly,
from $0.32$ at step 10 to $0.60$ at step 20 and $1.13$ at step 50.
(b) Test trajectory and EDMD
prediction in the first two DM coordinates; the rollout tracks the test
trajectory over the short horizon and drifts away over the long horizon.
The short-horizon closure is valid; the long-horizon error accumulates
system noise and dictionary truncation.}
\label{fig:edmd}
\end{figure}

A degree-2 EDMD operator is trained on the certified $n=30$ point and
rolled out for 50 steps from a held-out initial condition. The
dictionary $\phi(\psi)$ consists of the constant, the first two
diffusion-map coordinates $\psi_1,\psi_2$, and their quadratic terms
$\psi_1^2,\psi_1\psi_2,\psi_2^2$, giving six basis functions in total.
The EDMD matrix $A$ is fit by ridge regression with $\lambda=10^{-4}$ on
pairs $(\phi(\psi_t),\phi(\psi_{t+1}))$ from four training trajectories,
separate from the twelve trajectories used for the diffusion-map
embedding. The reported RMSE is the coordinate-wise root-mean-square
error on the two slow DM coordinates $\psi_1,\psi_2$, not normalized
against a trajectory scale; each rollout step corresponds to a stored
sample at $0.02$ time units, so 50 steps span about $1.0$ time unit.
The slow-coordinate RMSE grows
approximately linearly, from $0.32$ at step 10 to $0.60$ at step 20 and
$1.13$ at step 50 (Figure~\ref{fig:edmd}). The
drift follows from the setup: the certified two-dimensional manifold
supports the short-horizon closure, and long-horizon error accumulates
system noise and the truncation of the degree-2 dictionary. The
short-horizon closure holds precisely where the certificate places it.

\subsection{FitzHugh--Nagumo}

At the FHN Hopf point, the diffusion-map spectrum is
$[1.000,0.9947,0.9942,0.9800,\dots]$ at $n=30$, again showing two
near-degenerate leading components. At $n=192$, the spectrum is
$[1.000,0.9998,0.9991,0.9980,\dots]$: the leading components remain
near-degenerate but the gaps shrink to $10^{-3}$, consistent with the
certificate approaching its boundary.

Table~\ref{tab:fhn} summarizes the FHN spectrum across lattice sizes,
alongside the certificate ratio.

\begin{table}[t]
\caption{FitzHugh--Nagumo near Hopf: leading diffusion-map eigenvalues
and certificate ratio. The first DM eigenvalue is normalized to one and
omitted; the certified range extends to $n_*\approx 203$.}
\label{tab:fhn}
\begin{tabular}{r r r r r}
$n$ & $\lambda_1^{\rm DM}$ & $\lambda_2^{\rm DM}$ &
$\lambda_3^{\rm DM}$ & gap ratio\\
\hline
30  & 0.9947 & 0.9942 & 0.9800 & 2130\\
60  & 0.9973 & 0.9972 & 0.9894 & 134\\
100 & 0.9996 & 0.9983 & 0.9959 & 17.2\\
150 & 0.9990 & 0.9989 & 0.9961 & 3.37\\
192 & 0.9998 & 0.9991 & 0.9980 & 1.24\\
250 & 0.9998 & 0.9993 & 0.9984 & 0.42
\end{tabular}
\end{table}

The FHN certificate ratio decays from $2130$ at $n=30$ to $1.24$ at
$n=192$ and crosses $1$ near the closed-form estimate
$n_*^{\rm approx}\approx 203$. The DM eigenvalues
remain near-degenerate throughout the certified range, and the gap
between $\lambda_2^{\rm DM}$ and $\lambda_3^{\rm DM}$ declines from
$0.014$ at $n=30$ to $0.0011$ at $n=192$, a smaller dynamic
range than the Brusselator. The smaller FHN Lipschitz constant
($1.4$ vs $184$) opens a wider certified window, so the
spectral margin changes more gradually.

\section{Beyond the certified boundary}

Table~\ref{tab:beyond} extends the Brusselator spectrum to the diagnostic
arm $n=80,113$, where the certificate fails.

\begin{table}[t]
\caption{Brusselator beyond the certified boundary. The first DM
eigenvalue is normalized to one and omitted.}
\label{tab:beyond}
\begin{tabular}{r r r r r}
$n$ & $\lambda_1^{\rm DM}$ & $\lambda_2^{\rm DM}$ &
$\lambda_3^{\rm DM}$ & gap ratio\\
\hline
60  & 0.985 & 0.970 & 0.949 & 1.01\\
80  & 0.989 & 0.986 & 0.957 & 0.32\\
113 & 0.991 & 0.981 & 0.971 & 0.10
\end{tabular}
\end{table}

At $n=113$ the leading eigenvalues form a dense block
$0.991,0.981,0.971,0.955$ rather than a clean two-plus-one pattern. The
spectral-gap condition no longer holds; the data
spectrum shows the two-dimensional slow structure losing its clean
separation. The two readings agree.

Breaking the diffusion symmetry
($D_2=0.3,0.1$ instead of $1.0$) inside the certified range produces no
nontrivial transverse graph: the diffusion-map spectrum is
unchanged to four decimal places. On a ring with $n\le 60$, the smallest
transverse damping $\lambda_1\ge 27$ dwarfs the slow-mode growth
$\gamma=0.1$, so strong diffusion locks in the uniform limit cycle
regardless of the $D_1/D_2$ ratio. The trivial graph is a structural
feature of the certified regime, not a numerical artifact.

Table~\ref{tab:robust} collects the diagnostic-arm scan: nine
configurations $(n,D_2)$ show that $D_2$ leaves the leading DM eigenvalues
unchanged at any lattice size, while the leading modes tighten into a
dense block as $n$ grows
(Table~\ref{tab:beyond} lists the full leading spectrum).

\begin{table}[t]
\caption{Robustness scan: leading DM eigenvalue $\lambda_1^{\rm DM}$
across $(n,D_2)$; only $\lambda_1^{\rm DM}$ is listed, and the full
leading spectrum appears in Table~\ref{tab:beyond}. $D_2$ has no effect
at any $n$.}
\label{tab:robust}
\begin{tabular}{r r r r}
$n$ \textbackslash $D_2$ & 1.0 & 0.3 & 0.1\\
\hline
80  & 0.9893 & 0.9893 & 0.9893\\
113 & 0.9915 & 0.9915 & 0.9915\\
150 & 0.9945 & 0.9945 & 0.9945
\end{tabular}
\end{table}

\begin{figure}[t]
\centering
\includegraphics[width=0.8\columnwidth]{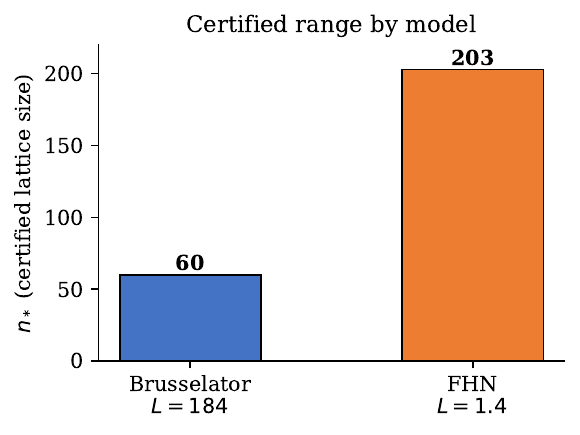}
\caption{The certified lattice size $n_*$ for the two model systems. The
reaction Lipschitz constant $L_{\Box}$ differs by two orders of magnitude
and fixes the certified range through Eq.~\eqref{eq:nstar}.}
\label{fig:nstar_models}
\end{figure}

\begin{figure}[t]
\centering
\includegraphics[width=0.9\columnwidth]{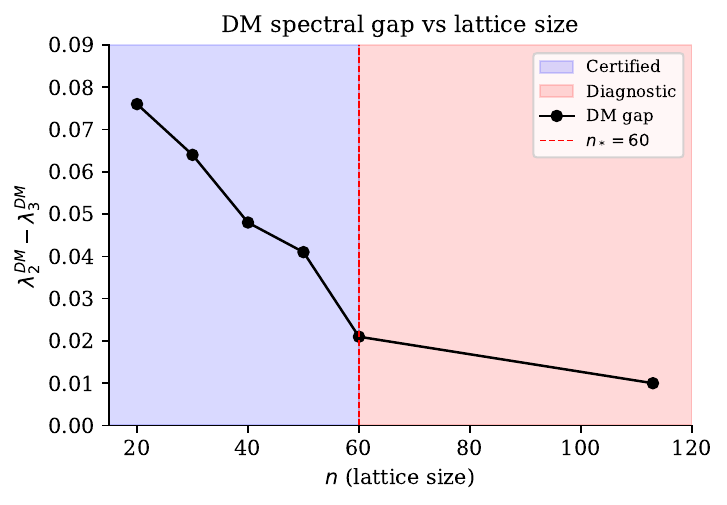}
\caption{The diffusion-map spectral gap $\lambda_2^{\rm DM}-\lambda_3^{\rm DM}$
as a function of lattice size $n$ on the Brusselator at $b=2.2$. The
shaded region is the certified range $n\le n_*\approx 60$; the gap shrinks
monotonically and the spectrum becomes dense beyond it.}
\label{fig:gap}
\end{figure}

\begin{figure*}[t]
\centering
\includegraphics[width=0.95\textwidth]{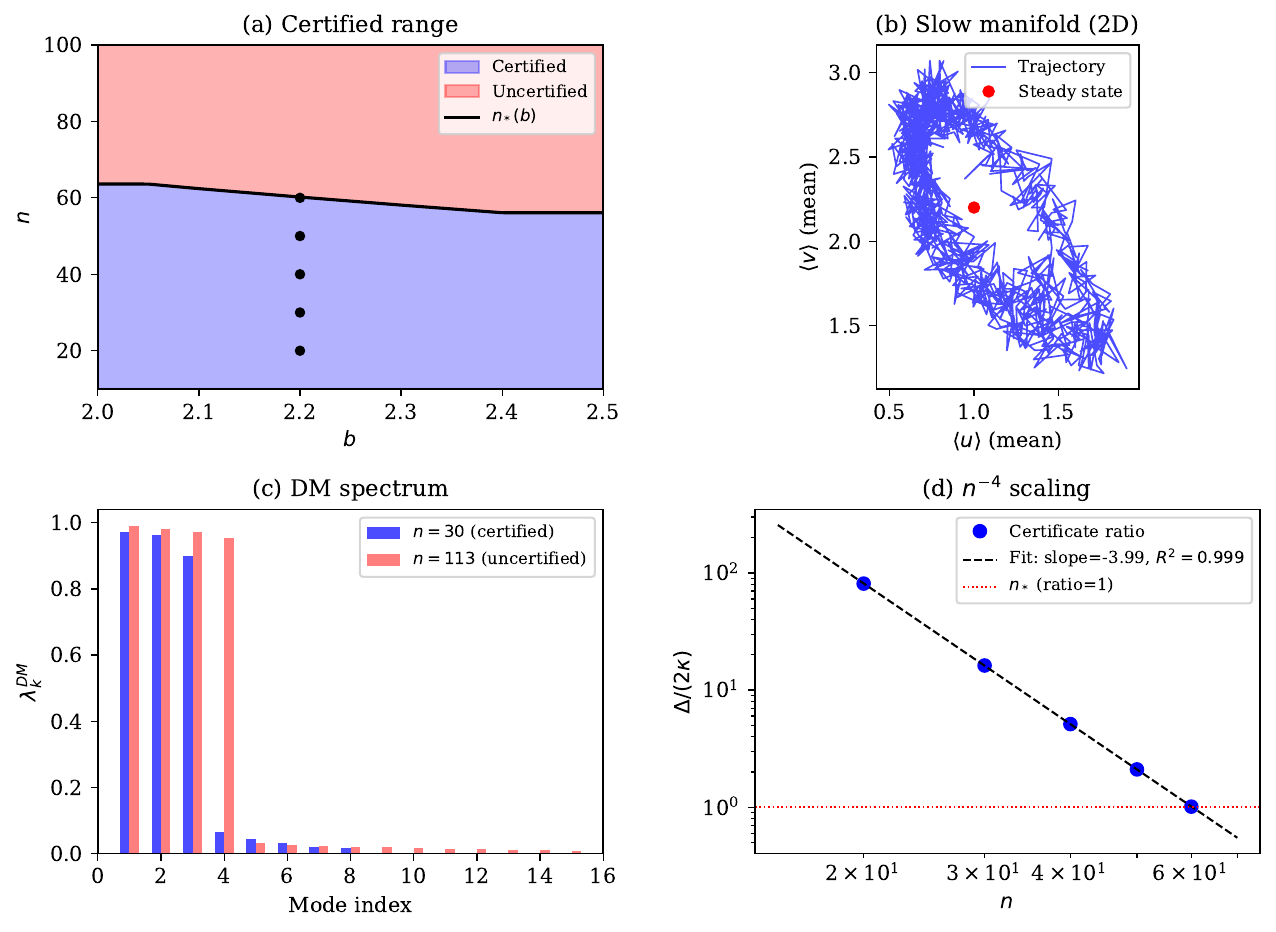}
\caption{Summary of the certification on the Brusselator lattice at
$b=2.2$. (a) Certified (blue) and uncertified (red) regions in the
$(n,b)$ plane; the solid curve is the numerical boundary $n_*(b)$.
Black dots mark the five certified operating points.
(b) Projection of the simulated trajectory onto the mean $(u,v)$
coordinates, showing the two-dimensional Hopf limit cycle (blue) and
the steady state (red dot). The trajectory is from a single certified
$n=30$ simulation at $b=2.2$.
(c) Diffusion-map spectrum at $n=30$ (certified, blue) and $n=113$
(uncertified, red). Below $n_*$ the spectrum shows a clean two-plus-one
separation; above $n_*$ the leading modes become near-degenerate and the
transverse modes densify.
(d) Log-log plot of the certificate margin $\Delta/(2\kappa)$ versus
$n$. The fitted slope is $-3.99\pm 0.08$ ($R^2=0.999$), consistent
with the predicted $n^{-4}$ scaling. The dotted line at ratio $=1$
marks $n_*$.}
\label{fig:main}
\end{figure*}

\begin{figure*}[t]
\centering
\includegraphics[width=0.95\textwidth]{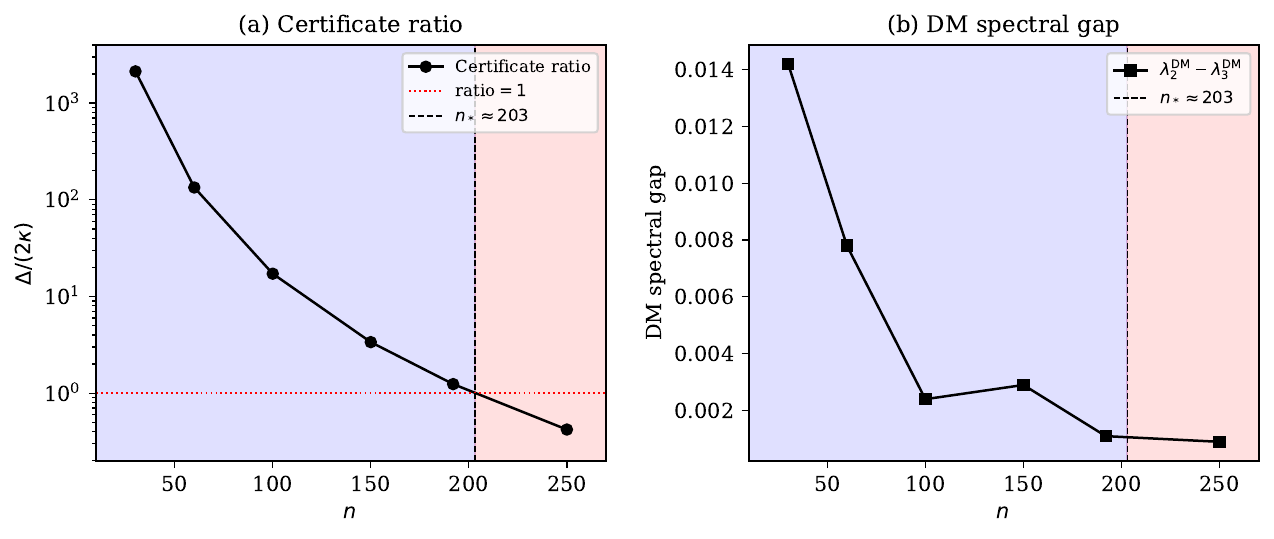}
\caption{FitzHugh--Nagumo near Hopf. (a) Certificate ratio
$\Delta/(2\kappa)$ versus lattice size $n$; the dotted line marks
ratio $=1$ and the dashed line the closed-form boundary
$n_*\approx203$. (b) DM spectral gap
$\lambda_2^{\rm DM}-\lambda_3^{\rm DM}$ on the same six points. The
certified range lies left of $n_*$; the theoretical margin and the DM
gap both decay toward the boundary.}
\label{fig:fhn}
\end{figure*}

\section{Discussion}

The certificate does not prove a nontrivial transverse graph exists; on the
Brusselator at $b=2.2$ the graph is the trivial embedding of the uniform
Hopf cycle, and this is reported directly. Its role is evidential: it
tells a practitioner running diffusion maps on a simulated lattice whether
the recovered two-dimensional object is the certified slow manifold or a
finite-resolution construction.

The certified range follows the reaction Lipschitz constant, not the
particular model. The FHN system, with $L_{\Box}$ two orders of
magnitude smaller than the Brusselator, reaches
$n\approx 203$; the Brusselator
stops at $n=60$. Eq.~\eqref{eq:nstar} estimates the certified resolution
before any learning algorithm runs.

The empirical spectrum and the certificate deteriorate together.
Below $n_*$ the leading diffusion-map eigenvalues show a clean two-plus-one
pattern; above $n_*$ that pattern loses its guarantee. In the
diagnostic arm the leading eigenvalues densify, though the spectrum need
not fill in at every size. The certificate predicts where the
two-dimensional reading is trustworthy.

The trivial-graph degeneracy is itself a result. On a strongly
damped Hopf lattice, the slow manifold is the cycle itself, and any data
method that recovers a two-dimensional structure recovers the cycle.
Claims of nontrivial inertial structure on such lattices should be read
against the certificate: above $n_*$, apparent extra dimensions
may be finite-resolution modes rather than certified slow directions.

The certificate is conservative. The clean-to-dense
transition sits near $n_*$ but not exactly at it, and the interval
cross-check at $n=60$ leaves only a 1\% margin. A tighter bound would
require a sharper estimate of the reaction coupling than the box
sup norm used here.

The certificate carries several restrictions. It covers periodic rings with equal
diffusion coefficients only: the maximum-principle chain defining the
absorbing box requires
$D_1=D_2$, and unequal diffusion or Neumann boundaries call for a
re-derived box invariant, a separate
calculation. It bounds manifold existence, not the accuracy of a
particular learning algorithm; a poor
dictionary or an unlucky bandwidth can still miss the
two-dimensional structure even inside $n_*$. It is
deterministic and does not account for stochastic noise in the dynamics
itself, only observation noise in the data; a stochastic version
would replace the spectral-gap condition with a mean-square or
almost-sure variant, outside the present scope. In
the Brusselator example the Hopf limit cycle is spatially uniform, so the
certified manifold is the cycle itself with trivial transverse graph
$\varphi\equiv0$; the two-dimensional reading largely recovers the periodic
phase rather than a nontrivial spatial slow mode. A near-two-dimensional
data spectrum and a valid existence guarantee are distinct:
the PCA energy and the correlation dimension stay near the cycle value
across the boundary, while only the DM gap tracks the certificate margin.
Data-side spectral degeneracy can persist after the guarantee lapses, and
a clean spectrum inside $n_*$ does not by itself prove that the manifold
exists.

Koopman convergence theory~\cite{Colbrook2026} characterizes in
the abstract when spectral learning can or cannot converge. Its framework
is system-agnostic; the present certificate instantiates one concrete case for
reaction-diffusion lattices, where the spectral-gap condition evaluates
from the reaction Lipschitz constant and returns a closed-form
resolution bound. The abstract theory marks which conditions to seek;
the certificate marks where those conditions hold on a concrete PDE
discretization. The earlier rigorous data-driven Koopman computation of
Colbrook and Townsend~\cite{ColbrookTownsend2024} sits at the same
system-agnostic level: ResDMD returns spectra and pseudospectra with
explicit convergence bounds and a posteriori error control, but its
bounds sit on the observable space rather than on the PDE discretization.
In its high-dimensional molecular and turbulent examples, a second
snapshot set checks a learned dictionary a posteriori for sufficiency;
that check still asks how well the dictionary resolves an existing
object, not whether the discretization supports one. The EDMD check here
keeps the same separation: the four trajectories fitting the operator
differ from the twelve used for the embedding, so the closure error
is not measured on the snapshots that selected the coordinates.
Colbrook \textit{et al.} match upper and lower bounds to separate
learnable from inadmissible regimes; this certificate supplies one
sufficient side only. It certifies existence below $n_*$, and a violation
above $n_*$ is a missed guarantee, not an impossibility
statement; a two-sided existence boundary on the lattice remains
open.

Conradie \textit{et al.}\cite{Conradie2026} develop
invariance diagnostics and error bounds for Koopman operator learning,
quantifying how a finite-dimensional feature dictionary fails to be
invariant. Their analysis takes the observable space as given: a
dictionary is fixed and the closure error is reported. The certificate here
takes the PDE side: at a given lattice resolution, it asks whether a
two-dimensional slow manifold exists at all. The two sides are both required.
A learned model is trustworthy only when the PDE side guarantees a
low-dimensional object exists and the observable side guarantees the
dictionary resolves it. The certificate supplies the PDE-side half of a
certification that the data-driven community has so far built from the
observable side alone.

Conradie \textit{et al.} divide the errors of
data-driven Koopman learning into estimation error from finite noisy
data, numerical error from finite-precision algebra, and projection error
from truncating the infinite-dimensional observable
space; their diagnostics control the second and third. Sampling-error
convergence rates for EDMD, gEDMD, and kernel EDMD are established under
i.i.d. or ergodic sampling~\cite{Philipp2024,Llamazares2026,Olguin2025};
these bounds likewise sit on the observable side, conditioning the
dictionary and the data rather than the lattice discretization. The certificate
here targets a failure mode that precedes all three. When
the discretization itself does not support a finite-dimensional slow
manifold, no dictionary, however well chosen, can recover the object,
because none exists at that resolution. This \emph{existence error}
is detectable before any learning runs,
from the reaction Lipschitz constant and the lattice spacing alone,
and it does not depend on the dictionary choice.

On the observable side, the recovered slow coordinates carry the three
errors listed above: estimation from the finite noisy snapshot set, numerical
from finite-precision ridge regression, and projection from the six-term
degree-2 dictionary. Their size depends on the sampling density and
ergodicity, the observation-noise level, the diffusion-map bandwidth, and the
dictionary degree; the EDMD rollout in Sec.~IV measures the
accumulated coordinate error empirically, not through a closed bound. The
certificate fixes only the PDE-side precondition, that a two-dimensional slow
object exists at the chosen $n$. It does not bound these observable-side
errors, and a small rollout RMSE inside $n_*$ does not tighten the existence
certificate in return.

Alla, Geelen, and Lu~\cite{Alla2026} report non-monotonic convergence
of POD-reduced models for
reaction-diffusion systems: adding reduced coordinates does not
monotonically lower the ROM error. This matches the
existence error identified here. Above $n_*$, the
additional POD modes resolve not a slow manifold but
finite-resolution transverse modes that grow with $n$. The certificate
places the onset of non-monotonicity near $n_*$, where the spectral gap
closes.

Data-driven continuation of patterns and their
bifurcations~\cite{DataContin2025} traces transition curves in
reaction-diffusion parameter space by Wasserstein distances between
pattern statistics. It locates pattern boundaries; the certificate here
answers the prior question of whether a low-dimensional object exists at
the chosen resolution. A continuation run above $n_*$ may
track finite-resolution artifacts rather than true pattern
boundaries: the spectral-gap condition has failed, and the
apparent low-dimensional structure carries no guarantee. The existence
check precedes any reliable bifurcation tracking.

Conditional generative models sample initial states on a slow
manifold and fill missing segments of a bifurcation
diagram~\cite{Crabtree2025}. The construction assumes the low-dimensional
object exists to be sampled; its conditioning label itself comes from a
manifold-learning step on the training data, the observable side the
certificate checks beforehand. Above $n_*$, the structure the
generative model conditions on may be a finite-resolution transverse
object rather than an inertial manifold.

Spectral submanifolds (SSMs) take a geometric route to the same
reduction problem. Liu, Ax{\aa}s, and Haller~\cite{HallerSSM2024} continue
a non-resonant spectral subspace at an equilibrium or periodic orbit
smoothly into a local $C^\infty$ graph, then fit a
polynomial reduced dynamics on that graph. Their examples range from
delay-embedded Lorenz and R\"ossler attractors to spatially extended
systems: a forced Duffing oscillator chain and the
Kuramoto--Sivashinsky equation. The SSM construction relies on a known
anchor state and a chosen slow eigenspace; its existence rests on the
local spectral separation at that anchor, not on a resolution bound
for the discrete lattice.

An SSM model reports how
well reduced dynamics reproduces short-time trajectories and
long-time statistics, such as the largest Lyapunov exponent and the
invariant distribution; this certificate reports whether a
finite-dimensional slow manifold exists on the given lattice. The two answers
differ in kind. They are compatible inside the certified range. Below $n_*$ the
lattice-level gap holds, and an SSM learned there attaches to a real
inertial object; above $n_*$, an SSM-like manifold fit on the lattice
tracks finite-resolution transverse modes rather than an inertial
manifold.

Xie and Rowley~\cite{XieRowley2026} embed the inertial-manifold prior
directly into a neural operator: IMNO uses the assumed
low-dimensional slow structure as an inductive bias for long-horizon
autoregressive prediction on dissipative PDEs. Their construction takes
the existence of a low-dimensional object as given and encodes it in the
network. The certificate supplies the precondition such encoding can
rely on. Below $n_*$ the lattice admits the inertial object the operator
assumes; above $n_*$, the same architecture fits
finite-resolution transverse modes. The certificate decides where an
inertial-manifold-based operator is meaningful.

Koronaki \textit{et al.}~\cite{Koronaki2024} organize data-driven ROMs
for dissipative PDEs around the approximate-inertial-manifold template.
On reaction-diffusion and Kuramoto--Sivashinsky systems they estimate
manifold dimension and latent variables from POD, diffusion maps, or
autoencoders, learn the approximate inertial form, and correct a
low-order Galerkin truncation by post-processing. The workflow assumes
an inertial manifold, or a thin AIM layer around the attractor, from the
outset; its data methods approximate and correct that object but do not
decide whether it exists at the chosen lattice resolution. The certificate
sits ahead of that workflow: it evaluates the spectral-gap condition on
the discrete lattice and returns $n_*$, below which the inertial object
the AIM construction relies on is certified to exist. Their results also
show that granting an AIM dimension does not by itself keep a truncated
low-order model accurate over the long term; when the post-processing
premise fails, a Gray-Box learned correction to the reduced dynamics
enters. The certificate addresses the question preceding that correction:
whether the discretization admits the low-dimensional object at all.

Conformal prediction has been combined with Koopman control,
giving distribution-free probabilistic error bounds on trajectory
prediction~\cite{ConformalKoop2026}. Those bounds are statistical and
hold at a user-chosen probability; this certificate is deterministic
and returns a hard resolution threshold. The two are compatible: a
conformal Koopman model certified below $n_*$ inherits the
PDE-side existence guarantee and the observable-side statistical
guarantee.

Loh and Mireles James~\cite{Loh2026} certify periodic solutions
recovered by a low-dimensional proxy against the infinite-dimensional
target, using a-posteriori computer-assisted proofs. Their certificate
runs after the proxy has been regressed: it asks whether a periodic
orbit found in the proxy has a nearby true orbit in the target. This
certificate runs before any regression: it asks whether the target
discretization admits a finite-dimensional slow manifold at all. The two
stages act in sequence. A proxy built above $n_*$ may pass or
fail the proxy-to-target proof, but the object it approximates carries
no inertial-manifold guarantee. The pre-hoc existence check
here tells the practitioner whether a post-hoc certification
of the Loh--Mireles James type is meaningful at the chosen resolution.
Their a-posteriori proof is likewise one-sided: a failed certification
does not rule out a true periodic orbit nearby, since the available bounds
may be insufficient, just as a violation above $n_*$ here is a missed
guarantee rather than an impossibility statement.

Diffusion maps and EDMD recover low-dimensional coordinates from
simulated reaction-diffusion trajectories. The certificate states the
conditions under which that recovery is trustworthy and decides whether
the two-dimensional coordinates track a true slow
manifold rather than a finite-resolution construction:
the reading holds below $n_*$, and $n_*$ computes
from the reaction parameters alone.

\section{Conclusion}

Data-driven dimension reduction on reaction-diffusion
lattices carries a computable sufficient boundary at a critical
size $n_*$, derived in closed form from the reaction Lipschitz
constant and the slow-mode growth rate. Below $n_*$, diffusion maps recover a
reading backed by the existence guarantee; above $n_*$, the guarantee
lapses, without ruling out a slow manifold. The
boundary recurs across nonlinear structure, and the
spectral margin decays as $n^{-4}$, set by the coupling between
slow and fast Fourier modes.

The classical inertial-manifold theory, developed decades ago for
purely mathematical purposes, supplies a prior that marks when to trust
the low-dimensional structure an algorithm discovers. On strongly damped
Hopf lattices the prior predicts the slow manifold collapses to the
uniform cycle itself, the dimensional degeneracy reported here.

For engineering practice, the certificate replaces two routine reduction
decisions with a bounded check. It
sets a resolution target: a simulation or sensor grid inside $n_*$ carries
the existence guarantee for the two-dimensional object the reduced
coordinates read off; above $n_*$ the certificate offers no such
guarantee, and a lapse of the guarantee does not establish that no slow
manifold exists. A reduced model is checked before use:
before trusting a POD, autoencoder, or Koopman surrogate on a
lattice, the lattice-level gap can be evaluated as an existence check at
the chosen resolution; a richer dictionary or an extra reduced coordinate
cannot substitute for that check. The
check tolerates the observation-noise levels tested here, but it assumes
equal diffusivities, a periodic ring, and deterministic dynamics, and
carries no guarantee beyond those assumptions. The gap and Lipschitz
inputs do not depend on the Brusselator or FHN reaction forms, so the same
check can act as a prior for other reaction-diffusion field systems, such
as catalytic surfaces, ecological patches, or neural fields, provided their
lattice meets the same diffusivity and boundary assumptions; this is a
precondition for reduction, not a certified case study.

Several directions remain open. Extending the transition from periodic
rings to two-dimensional lattices and unstructured meshes widens
its applicability. Coupling the PDE-side resolution guarantee
with the observable-space guarantee of abstract Koopman
convergence theory~\cite{Colbrook2026} would yield a two-sided
certification. The $n^{-4}$ scaling, derived here for
reaction-diffusion lattices, may generalize to other multiscale
systems where a slow manifold competes with a fast spectral gap.

A further direction extends the certificate to stochastic lattices.
The present certificate
is deterministic; observation noise perturbs the embedding but leaves
the manifold untouched. In systems with dynamical noise, additive or
multiplicative white noise in the lattice reaction terms, the
deterministic spectral-gap condition gives way to a mean-square or
almost-sure variant, and the inertial manifold becomes a random
inertial manifold~\cite{Bensoussan1995}. The computable framework here,
evaluating gap, box Lipschitz, and slow growth on a concrete
lattice, carries over directly: the gap and Lipschitz inputs stay
unchanged, and the growth rate $\gamma$ acquires a noise-induced shift
estimable from the noise intensity. Such an extension reaches the noisy
lattice systems that form the setting of most data-driven
simulations.

\section*{Author contributions}
Y.L. conceived the project, developed the certificate framework,
performed the numerical computations, and wrote the manuscript.

\section*{Competing interests}
The author declares no competing financial or non-financial interests.

\section*{Data availability}
The simulation and certificate codes, along with the JSON data files
supporting the reported findings, are available from the author
upon reasonable request.

\end{document}